# GLI ANGOLI ALLA BASE DI UN TRIANGOLO ISOSCELE

Negli ultimi mesi, anche parlando dell'esame di Stato a conclusione del Liceo Scientifico, si è molto insistito sul fatto che le Indicazioni Nazionali raccomandano di presentare le applicazioni della matematica e, in particolare, il concetto di modello matematico. Non vorrei però che si dimenticasse che le stesse Indicazioni, come pure le Linee Guida per gli Istituti Tecnici e Professionali, sottolineano d'altra parte l'importanza del metodo dimostrativo: «*Verrà chiarita l'importanza e il significato dei concetti di postulato, assioma, definizione, teorema, dimostrazione*» (per inciso, io non farei differenza fra postulato e assioma).

In questo breve articolo vorrei confrontare varie dimostrazioni di uno dei primi teoremi che si incontrano in geometria euclidea del piano: «*gli angoli alla base di un triangolo isoscele sono uguali*», ovvero «*se un triangolo ha due lati uguali (congruenti), allora gli angoli opposti a quei lati sono uguali (congruenti)*». Di questo teorema si è già parlato su Archimede vari anni fa ([2] e [3]).

Il teorema compare all'inizio degli *Elementi* di Euclide (libro I, proposizione 5), dove è immediatamente seguito dal suo inverso (proposizione 6). Forse anche per la sua collocazione all'inizio del curriculum di geometria, il teorema è talvolta indicato con il nome di *pons asinorum*.

L'enunciato è di immediata verifica: basta ritagliare un triangolo isoscele disegnato su un foglio di carta e ripiegarlo in modo da far combaciare i lati uguali. Una tale verifica, sicuramente adatta all'inizio di una Scuola secondaria di I grado, non è una dimostrazione, ma è comunque didatticamente più ricca di un semplice controllo dell'uguaglianza dei due angoli effettuato con un goniometro: la piegatura della carta fa capire che il triangolo isoscele ammette un asse di simmetria (anche se non si usa esplicitamente questa parola), e che segmenti o angoli che si corrispondono in una simmetria sono uguali.

Chi dimenticasse la collocazione iniziale del nostro teorema potrebbe facilmente proporre dimostrazioni che sembrano rapide e brillanti. Vediamo due esempi.

- Siano $a$, $a'$ le lunghezze di due lati di un triangolo e sia $b$ quella del terzo lato; indichiamo con $\alpha$, $\alpha'$ le ampiezze degli angoli opposti ad $a$, $a'$. L'area del triangolo si può esprimere sia come $(1/2)a'b \operatorname{sen} \alpha$ sia come $(1/2)ab \operatorname{sen} \alpha'$. Dall'uguaglianza $(1/2)a'b \operatorname{sen} \alpha = (1/2)ab \operatorname{sen} \alpha'$ segue subito che $a = a'$ se e solo se $\operatorname{sen} \alpha = \operatorname{sen} \alpha'$; e, tenendo presente che due angoli di uno stesso triangolo non possono essere supplementari, concludiamo $a = a'$ se e solo se $\alpha = \alpha'$, ottenendo così anche l'enunciato inverso.

- Dato un triangolo *ABC* con *AB* = *AC*, tracciamo la circonferenza di centro *A* e raggio *AB*, che passa anche per *C* (fig. 1). Gli angoli del triangolo in *B* e in *C* sono angoli alla circonferenza a cui corrispondono archi uguali (una semicirconferenza diminuita dell'arco *BC*): quindi i due angoli sono uguali.

*Tuttavia*, dimostrazioni del tipo precedente presentano due grossi difetti. In primo luogo bisogna controllare con cura che tutti i teoremi che si danno per noti si riescano a dimostrare senza presupporre il teorema che stiamo considerando. Per esempio, nelle dimostrazioni classiche dei teoremi relativi agli angoli al centro e alla circonferenza si usa proprio il teorema sugli angoli alla base di un triangolo isoscele.

Inoltre, il teorema di cui stiamo parlando vale anche nelle geometrie ellittica e iperbolica: se applichiamo teoremi specifici di geometria euclidea, troviamo un risultato che dipende dal postulato delle parallele. Così, il ricorso alla trigonometria o all'area di un triangolo (calcolata nel modo usuale), o il fatto di escludere che un triangolo abbia due angoli supplementari, limita il risultato all'ambito della geometria euclidea.

Esaminiamo ora tre dimostrazioni che si possono effettivamente presentare all'inizio della trattazione della geometria. In ogni caso, si presuppone il I criterio di uguaglianza dei triangoli (due lati e l'angolo compreso).

Sia, come prima, *ABC* un triangolo isoscele con i lati *AB* e *AC* uguali.

**Prima dimostrazione** – Si traccia la bisettrice *AH* e si considerano i triangoli *ABH* ed *ACH*. Questi due triangoli sono uguali per il I criterio; si conclude subito l'uguaglianza degli angoli alla base.

La dimostrazione è rapida e convincente, ma presenta due problemi. In primo luogo, a rigore, si deve dimostrare a priori l'esistenza della bisettrice. La costruzione della bisettrice di un angolo si trova nel libro I degli *Elementi* solo come proposizione 9: Euclide la ottiene applicando il III criterio di uguaglianza, che, a sua volta, dipende proprio dal teorema di cui ci stiamo occupando. Inoltre, sul piano didattico, non credo ci sia un procedimento altrettanto semplice per dimostrare il teorema inverso.

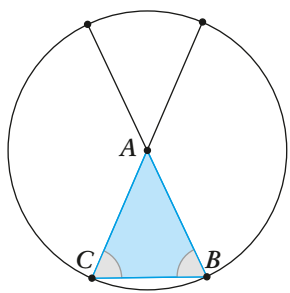

Figura **1**

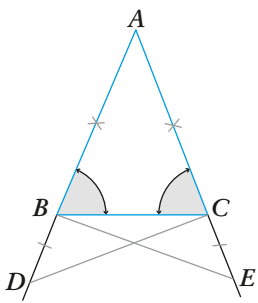

Figura **2**

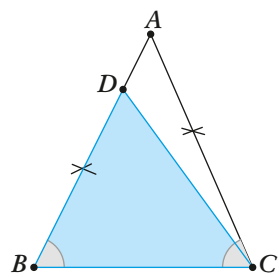

Figura **3**

**Seconda dimostrazione** – È la classica dimostrazione di Euclide. Si prolungano i lati uguali con due segmenti uguali *BD* e *CE*, scelti a piacere (fig. 2). Congiungendo *B* con *E* e *C* con *D*, si ottengono i triangoli *ACD* e *ABE* uguali per il I criterio di uguaglianza. Poi si considerano i triangoli *DBC* ed *ECB*, che risultano a loro volta si uguali, ancora per il I criterio. Se ne deduce che sono uguali gli angoli in *B* e in *C* di questi ultimi triangoli; e, dunque, anche gli angoli in *B* e in *C* del triangolo *ABC*.

Con la stessa costruzione, ma facendo riferimento sia al I sia al II criterio di uguaglianza, si dimostra il teorema inverso. Si considerano prima i triangoli *BDC* ed *ECB*, che risultano uguali per il I criterio; dall'uguaglianza degli altri elementi di questi due triangoli si deduce l'uguaglianza dei triangoli *ACD* e *ABE* per il II criterio, e si conclude che *AB* è uguale ad *AC*.

Aggiungo due osservazioni. In primo luogo, si può obiettare che queste dimostrazioni, per il teorema diretto e per il teorema inverso, sono troppo complesse per risultati così semplici e intuitivi. A mio parere, tuttavia, una dimostrazione non serve solo per spiegare e convincere, ma anche per sistemare e organizzare i risultati (si vedano [1] e [4]). Sarebbe molto pericoloso dare per scontati tutti i fatti intuitivamente ovvi (come indubbiamente è il nostro enunciato); inoltre, una dimostrazione è importante anche perché aiuta a capire le proprietà coinvolte e a porre il discorso nel giusto contesto.

In secondo luogo, notiamo che la dimostrazione di Euclide del teorema inverso è completamente diversa. Se per assurdo il lato *AB* fosse ad esempio maggiore del lato *AC*, allora Euclide considera il punto *D* di *AB* tale che *AD* sia uguale ad *AC* (fig. 3). I triangoli *ABC* e *DBC* avrebbero: *BC* in comune, *AC* = *DB*, *ACB* = *DBC*; dunque i due triangoli sarebbero uguali per il I criterio, ma questo è assurdo perché «il tutto è maggiore della parte», come recita la più celebre delle *nozioni comuni* di Euclide.

**Terza dimostrazione** – È la dimostrazione più rapida (*non* la più semplice), che pare risalga addirittura a Pappo di Alessandria (IV sec. d.C.). Si tratta del ragionamento precedente con *B* = *D* e *C* = *E*: cioè si considerano i due «triangoli» *ABC* e *ACB*, in quest'ordine.

Per essere rigorosi, occorre enunciare con maggiore cura il primo criterio di uguaglianza: date due *terne ordinate* di punti, *A*, *B*, *C* e *A'*, *B'*, *C'*, se valgono le uguaglianze *AB* = *A'B'*, *CA* = *C'A'*, *A* = *A'*, allora si ha: *BC* = *B'C'*, *B* = *B'*, *C* = *C'* (cioè sono ordinatamente uguali anche gli altri elementi dei due «triangoli»).

Se applichiamo il criterio in questa forma alle due terne *A*, *B*, *C* e *A*, *C*, *B*, arriviamo subito alla conclusione. Ho enunciato il criterio facendo riferimento alle terne ordinate di punti e non ai triangoli perché i triangoli *ABC* e *ACB* sono banalmente sempre uguali, perché sono lo *stesso* triangolo.

In modo del tutto analogo si dimostra il teorema inverso mediante il II criterio, a sua volta enunciato con riferimento a due terne ordinate di punti anziché a due triangoli.

Penso sia quasi superfluo osservare che questa terza dimostrazione è indubbiamente più rapida ed elegante delle precedenti, ma, altrettanto indubbiamente, non è adatta per la maggioranza degli studenti dei primi anni delle Superiori.

Infine, aggiungo che il precedente enunciato del I criterio richiama da vicino l'ultimo assioma di congruenza dei *Fondamenti della geometria* di Hilbert (1899). Nei primi assiomi di congruenza Hilbert postula semplici proprietà relative alla possibilità di trasportare e sommare segmenti e angoli; dopo di che conclude gli assiomi di congruenza richiedendo che: «se valgono le uguaglianze $AB = A'B'$, $CA = C'A'$, $A = A'$, allora si ha anche $B = B'$, $C = C'$». Sulla base di quest'ultimo assioma, Hilbert dimostra poi i criteri di congruenza dei triangoli nella forma consueta.

**Claudio Bernardi**
claudio.bernardi@uniroma1.it